# The Lie Representation of the Complex Unit Sphere

Jason Glowney
*Department of Applied Mathematics*
*University of Washington*

**Abstract**

We present the derivation of the 6-dimensional Eulerian Lie group of the form SO(3,C). We describe our derivation process, which involves the creation of a finite group by using permutation matrices, and the exponentiation of the adjoint form of the subset representing the generators of the finite group. We take clues from the 2-dimensional complex rotation matrix to present, what we believe, is a true representation of the Lie group for the six-dimensional complex unit sphere and proceed to study its dynamics. With this approach, we can proceed to present this SO(3,C) group and derive its unitary counterpart that is U(3). The following findings can prove useful in mathematical physics, complex analysis and applications in deriving higher dimensional forms of similar division algebras.

**Introduction**

We start with introducing permutations that represent the two-dimensional complex form of $e^{ix}$. This is an order 2 finite group and the following permutation matrices are the elements of the group

$$a = \begin{bmatrix} 1 & 0 \\ 0 & 1 \end{bmatrix} \text{ and } b = \begin{bmatrix} 0 & -1 \\ 1 & 0 \end{bmatrix}. \tag{1}$$

We can create the adjoint representation then as follows

$$\begin{bmatrix} a & -b \\ b & a \end{bmatrix}. \tag{2}$$

The permutation "a", represents the identity of the group, as it runs the diagonal. The "b" permutation represents the imaginary axis of the rotation matrix, as evidenced by when taking its exponent, we find it equal to the negative identity

$$\begin{bmatrix} 0 & -1 \\ 1 & 0 \end{bmatrix}^2 = \begin{bmatrix} -1 & 0 \\ 0 & -1 \end{bmatrix}. \tag{3}$$

The above is equivalent to $i^2 = -1$ and represents an important way to determine the nature of permutations when looking at the derivation of division algebras. We can go on to set up our permutations matrices and exponentiate separately and find

$$\exp\left(\begin{bmatrix} a & 0 \\ 0 & a \end{bmatrix} + \begin{bmatrix} 0 & -b \\ b & 0 \end{bmatrix}\right) = \begin{bmatrix} e^a & 0 \\ 0 & e^a \end{bmatrix}\begin{bmatrix} cosb & -sinb \\ sinb & cosb \end{bmatrix}. \tag{4}$$

From the above we have the rotation matrix for our two-dimensional complex equation $|e^{ix}| = 1$ in the form of

$$\begin{bmatrix} cosb & -sinb \\ sinb & cosb \end{bmatrix}. \tag{5}$$

We define the distance function of *(5)* as the determinant of the matrix. In this two-dimensional complex form ($C_2$), we find this distance function to be equal to unity, as we are describing the complex unit circle

$$\det\left(\begin{bmatrix} cosb & -sinb \\ sinb & cosb \end{bmatrix}\right) = \cos^2 b + \sin^2 b = 1. \tag{6}$$

## Derivation of the Complex Unit Sphere

With the aforementioned permutations of $C_2$ in mind, our goal is to create a group of permutation matrices, that once exponentiated will yield a 6x6 rotation matrix, with three imaginary permutations, spaced symmetrically, and with a determinate equal to unity. Typically, in an order 6 finite group, there are a total of 6! permutation matrices to choose from, representing 720 different possibilities.

The permutations we present below represent imaginary axes, as when squared they are equal to the negative identity

$$b = \begin{bmatrix} 0 & -1 & 0 & 0 & 0 & 0 \\ 1 & 0 & 0 & 0 & 0 & 0 \\ 0 & 0 & 0 & 0 & 0 & -1 \\ 0 & 0 & 0 & 0 & 1 & 0 \\ 0 & 0 & 0 & -1 & 0 & 0 \\ 0 & 0 & 1 & 0 & 0 & 0 \end{bmatrix} \rightarrow b^2 = \begin{bmatrix} -1 & 0 & 0 & 0 & 0 & 0 \\ 0 & -1 & 0 & 0 & 0 & 0 \\ 0 & 0 & -1 & 0 & 0 & 0 \\ 0 & 0 & 0 & -1 & 0 & 0 \\ 0 & 0 & 0 & 0 & -1 & 0 \\ 0 & 0 & 0 & 0 & 0 & -1 \end{bmatrix},$$

$$d = \begin{bmatrix} 0 & 0 & 0 & 0 & 0 & -1 \\ 0 & 0 & 0 & 0 & 1 & 0 \\ 0 & 0 & 0 & -1 & 0 & 0 \\ 0 & 0 & 1 & 0 & 0 & 0 \\ 0 & -1 & 0 & 0 & 0 & 0 \\ 1 & 0 & 0 & 0 & 0 & 0 \end{bmatrix} \rightarrow d^2 = \begin{bmatrix} -1 & 0 & 0 & 0 & 0 & 0 \\ 0 & -1 & 0 & 0 & 0 & 0 \\ 0 & 0 & -1 & 0 & 0 & 0 \\ 0 & 0 & 0 & -1 & 0 & 0 \\ 0 & 0 & 0 & 0 & -1 & 0 \\ 0 & 0 & 0 & 0 & 0 & -1 \end{bmatrix}, \qquad (7)$$

$$f = \begin{bmatrix} 0 & 0 & 0 & -1 & 0 & 0 \\ 0 & 0 & 1 & 0 & 0 & 0 \\ 0 & -1 & 0 & 0 & 0 & 0 \\ 1 & 0 & 0 & 0 & 0 & 0 \\ 0 & 0 & 0 & 0 & 0 & -1 \\ 0 & 0 & 0 & 0 & 1 & 0 \end{bmatrix} \rightarrow f^2 = \begin{bmatrix} -1 & 0 & 0 & 0 & 0 & 0 \\ 0 & -1 & 0 & 0 & 0 & 0 \\ 0 & 0 & -1 & 0 & 0 & 0 \\ 0 & 0 & 0 & -1 & 0 & 0 \\ 0 & 0 & 0 & 0 & -1 & 0 \\ 0 & 0 & 0 & 0 & 0 & -1 \end{bmatrix}.$$

We continue by setting the "a" permutation to be our identity and "c" and "e" to represent the last 2 permutation matrices to arrive at our adjoint representation $\mathfrak{G}$.

$$\mathfrak{G} = \begin{bmatrix} a & -b & c & -f & e & -d \\ b & a & f & c & d & e \\ e & -f & a & -d & c & -b \\ f & e & d & a & b & c \\ c & -d & e & -b & a & -f \\ d & c & b & e & f & a \end{bmatrix}. \qquad (8)$$

We continue by looking at the inter-relations of our permutations from above. We find that the "c" and "e" permutations represent real-valued elements, where taking the product of any two of our imaginary permutations, will yield either "-c" or "-e". This is akin to the imaginary values yielding "-a" when squared. It will turn out that the permutations we consider in "-c" and "-e", are important in that they represent left and right-handed rotations from the product of axes, and they can also represent scaling factors that can create a gradient for progressing rotations. It is easy verify that $\mathfrak{G}$ is closed under multiplication as we will soon demonstrate.

|   | a | b  | c  | d  | e  | f  |
|---|---|----|----|----|----|----|
| a | a | b  | c  | d  | e  | f  |
| b | b | -a | f  | -e | d  | -c |
| c | c | d  | e  | f  | a  | b  |
| d | d | -c | b  | -a | f  | -e |
| e | e | f  | a  | b  | c  | d  |
| f | f | -e | d  | -c | b  | -a |

**Fig. 1** – Standard Form Cayley Table of Permutation matrix *(8)*. The identity element is noted as "a" above.

We establish the relationships between the elements of the permutation matrix by way of a Standard Cayley table [1] seen in **Fig. 1**. The table is read by taking the respective row and multiplying by the column. So, as an example taking *db* will yield *-c* as shown above. The reader will notice the symmetry between the alternating signs of the values for "a", "b" and "c" as one moves along the columns from left to right. We note that by way of our Standard Form Cayley table, that this is a non-abelian form under multiplication, as it is lacking symmetry along the diagonal and that some products do not commute, as well as our permutations are closed under multiplication.

With quaternions and octonion forms, the general approach has been to assume there is one real value and three imaginary values for quaternions, and one real and seven imaginary values for octonions. To have the symmetry that we have in the 2-dimensional complex form, we require equal real parts to imaginary parts for our complex unit sphere derivation. The key in arriving at this form, is the realization that permutations "a", "c" and "e" represent real values that are oriented plane(s). This is akin to a Clifford algebra form where $e_1 e_2 = -e_2 e_1$. Instead we have the following relations

$$b^2 = d^2 = f^2 = -1$$

$$bd = -e, df = -e, fb = -e$$

$$bf = -c, fd = -c, db = -c \qquad (9)$$

$$bf = fd, df = fb, bd = df$$

$$e^3 = -1, c^3 = -1, ce = 1, (-c)(-e) = 1.$$

We can use the relations above and simplify our algebras. For example

$$bfddbdfb \to -bfbdfb \to -bfbddf \to bfbf \to bffd \to -bd = e. \qquad (10)$$

So, the "c" and "e", are these oriented plane(s). In our permutation matrix, the imaginary permutations "b", "d" and "f" comprise the basis of the group 𝔊, by way of the fact that with the products of just these three permutations, we can arrive at all of the permutations in our finite group. We cannot do the same with the "a", "c" and "e" permutations. This is similar to using only the products between 2 imaginary values, you can arrive at any real and or imaginary number, but we cannot take real numbers and arrive at imaginary ones with the product of integers only. For example, the real number "4" can be obtained by taking the product of -2i and 2i. With these ideas in mind, we set to zero the permutations "a", "c" and "e" from the matrix *(8)* above and we denote it as 𝔤. (Note that this is similar to setting "a" = 0 in the quaternion form)

$$\mathfrak{g} = \begin{bmatrix} 0 & -b & 0 & -f & 0 & -d \\ b & 0 & f & 0 & d & 0 \\ 0 & -f & 0 & -d & 0 & -b \\ f & 0 & d & 0 & b & 0 \\ 0 & -d & 0 & -b & 0 & -f \\ d & 0 & b & 0 & f & 0 \end{bmatrix}. \qquad (11)$$

In the parlance of group theory, we state that g is a subset of 𝔊 ($g \in \mathfrak{G}$).

## SO(3,C) Lie Group and Algebras

The commutators of the elements in our associative algebra, where $[a,b] = ab - ba$, are as follows

$$[b,d] = c - e, [d,b] = e - c, [d,f] = c - e, [f,d] = e - c, [f,b] = c - e, [b,f] = e - c,$$

$$[b,c] = f - d, [c,b] = d - f, [d,c] = b - f, [c,d] = f - b, [f,c] = d - b, [c,f] = b - d, \quad (12)$$

$$[b,e] = d - f, [e,b] = f - d, [d,e] = f - b, [e,d] = b - f, [f,e] = b - d, [e,f] = d - b,$$

$$[c,e] = 0, [e,c] = 0.$$

To aid in obtaining a familiar feel, we substituted $x = b, y = d$ and $z = f$ into the above to arrive at

$$\mathfrak{g} = \begin{bmatrix} 0 & x & 0 & z & 0 & y \\ -x & 0 & -z & 0 & -y & 0 \\ 0 & z & 0 & y & 0 & x \\ -z & 0 & -y & 0 & -x & 0 \\ 0 & y & 0 & x & 0 & z \\ -y & 0 & -x & 0 & -z & 0 \end{bmatrix}. \quad (13)$$

We can verify that the permutations that form our adjoint matrix 𝔊, form a Lie algebra by way of [2]

- Bilinearity: $[aX + bY, Z] = a[X,Z] + b[Y,Z]$ and $[Z, aX + bY] = a[Z,X] + b[Z,Y]$ for $\forall \, X,Y,Z \in \mathfrak{g}$ and arbitrary numbers $a$ & $b$,

- Anticommutatvity: $[X,Y] = -[Y,X], [X,Z] = -[Z,X], [Y,X] = -[X,Y], [Y,Z] = -[Z,Y], [Z,X] = -[X,Z], [Z,Y] = -[Y,Z] \; \forall \, x,y,z \in \mathfrak{g}$,

- Jacobi Identity: $[X,[Y,Z]] + [Z,[X,Y]] + [Y,[Z,X]] = 0 \; \forall \, X,Y,Z \in \mathfrak{g}$.

We find that our matrices X, Y and Z are skew-symmetric where

$$X^T + X = 0, Y^T + Y = 0, Z^T + Z = 0 \text{ and } (X + Y + Z)^T + (X + Y + Z) = 0. \quad (14)$$

For example

$$\mathfrak{g}^T + \mathfrak{g} = \begin{bmatrix} 0 & x & 0 & z & 0 & y \\ -x & 0 & -z & 0 & -y & 0 \\ 0 & z & 0 & y & 0 & x \\ -z & 0 & -y & 0 & -x & 0 \\ 0 & y & 0 & x & 0 & z \\ -y & 0 & -x & 0 & -z & 0 \end{bmatrix} + \begin{bmatrix} 0 & -x & 0 & -z & 0 & -y \\ x & 0 & z & 0 & y & 0 \\ 0 & -z & 0 & -y & 0 & -x \\ z & 0 & y & 0 & x & 0 \\ 0 & -y & 0 & -x & 0 & -z \\ y & 0 & x & 0 & z & 0 \end{bmatrix} = 0. \quad (15)$$

As our matrix g is skew-symmetric, we are assured that with exponentiation, we will arrive at an orthogonal Lie group representation as

$$G = \exp(\mathfrak{g}) = \sum_{n=0}^{\infty} \frac{\mathfrak{g}^n}{n!}. \quad (16)$$

$$\begin{bmatrix} \frac{2\cos(r)}{3}+\frac{\cos(\gamma)}{3} & \frac{\partial r}{\partial x}\cdot\frac{2\sin(r)}{3}+\frac{\sin(\gamma)}{3} & -\frac{\cos(r)}{3}+\frac{\cos(\gamma)}{3} & \frac{\partial r}{\partial z}\cdot\frac{2\sin(r)}{3}+\frac{\sin(\gamma)}{3} & -\frac{\cos(r)}{3}+\frac{\cos(\gamma)}{3} & \frac{\partial r}{\partial y}\cdot\frac{2\sin(r)}{3}+\frac{\sin(\gamma)}{3} \\ -\frac{\partial r}{\partial x}\cdot\frac{2\sin(r)}{3}-\frac{\sin(\gamma)}{3} & \frac{2\cos(r)}{3}+\frac{\cos(\gamma)}{3} & -\frac{\partial r}{\partial z}\cdot\frac{2\sin(r)}{3}-\frac{\sin(\gamma)}{3} & -\frac{\cos(r)}{3}+\frac{\cos(\gamma)}{3} & -\frac{\partial r}{\partial y}\cdot\frac{2\sin(r)}{3}-\frac{\sin(\gamma)}{3} & -\frac{\cos(r)}{3}+\frac{\cos(\gamma)}{3} \\ -\frac{\cos(r)}{3}+\frac{\cos(\gamma)}{3} & \frac{\partial r}{\partial z}\cdot\frac{2\sin(r)}{3}+\frac{\sin(\gamma)}{3} & \frac{2\cos(r)}{3}+\frac{\cos(\gamma)}{3} & \frac{\partial r}{\partial y}\cdot\frac{2\sin(r)}{3}+\frac{\sin(\gamma)}{3} & -\frac{\cos(r)}{3}+\frac{\cos(\gamma)}{3} & \frac{\partial r}{\partial x}\cdot\frac{2\sin(r)}{3}+\frac{\sin(\gamma)}{3} \\ -\frac{\partial r}{\partial z}\cdot\frac{2\sin(r)}{3}-\frac{\sin(\gamma)}{3} & -\frac{\cos(r)}{3}+\frac{\cos(\gamma)}{3} & -\frac{\partial r}{\partial y}\cdot\frac{2\sin(r)}{3}-\frac{\sin(\gamma)}{3} & \frac{2\cos(r)}{3}+\frac{\cos(\gamma)}{3} & -\frac{\partial r}{\partial x}\cdot\frac{2\sin(r)}{3}-\frac{\sin(\gamma)}{3} & -\frac{\cos(r)}{3}+\frac{\cos(\gamma)}{3} \\ -\frac{\cos(r)}{3}+\frac{\cos(\gamma)}{3} & \frac{\partial r}{\partial y}\cdot\frac{2\sin(r)}{3}+\frac{\sin(\gamma)}{3} & -\frac{\cos(r)}{3}+\frac{\cos(\gamma)}{3} & \frac{\partial r}{\partial x}\cdot\frac{2\sin(r)}{3}+\frac{\sin(\gamma)}{3} & \frac{2\cos(r)}{3}+\frac{\cos(\gamma)}{3} & \frac{\partial r}{\partial z}\cdot\frac{2\sin(r)}{3}-\frac{\sin(\gamma)}{3} \\ -\frac{\partial r}{\partial y}\cdot\frac{2\sin(r)}{3}-\frac{\sin(\gamma)}{3} & -\frac{\cos(r)}{3}+\frac{\cos(\gamma)}{3} & -\frac{\partial r}{\partial x}\cdot\frac{2\sin(r)}{3}-\frac{\sin(\gamma)}{3} & -\frac{\cos(r)}{3}+\frac{\cos(\gamma)}{3} & -\frac{\partial r}{\partial z}\cdot\frac{2\sin(r)}{3}+\frac{\sin(\gamma)}{3} & \frac{2\cos(r)}{3}+\frac{\cos(\gamma)}{3} \end{bmatrix}$$

***Fig. 2*** – Lie Group Matrix ***G***. $r = \sqrt{x^2 - xy - xz + y^2 - yz + z^2}$, $\gamma = (x + y + z)$, $\frac{\partial r}{\partial x} = \frac{(x - \frac{y}{2} - \frac{z}{2})}{r}$, $\frac{\partial r}{\partial y} = \frac{(-\frac{x}{2} + y - \frac{z}{2})}{r}$, $\frac{\partial r}{\partial z} = \frac{(-\frac{x}{2} - \frac{y}{2} + z)}{r}$

We can now proceed to exponentiate the matrix g *(13)* to arrive at the Lie group that represents $|e^{ix+iy+iz}| = 1$. We utilized MATLAB and the *expm(*g*)* function to determine the result and after considerable manual simplification, a consistent pattern emerged in the form of our final result seen in ***Fig. 2.***

We see our rotation matrix is a mixture of two different arguments that the trigonometric functions operate on in $r = \sqrt{x^2 - xy - xz + y^2 - yz + z^2}$ and $\gamma = (x + y + z)$. With this rotation group, we have effectively split the exponential $e^{i(x+y+z)}$ into six dimensions, with 3 real and 3 imaginary axes. All of the columns of the matrix above are orthonormal to one another. We have in *'r'*, a non-Euclidian radius that acts in tandem with the $'\gamma'$ argument of the other trigonometric entries. The dual arguments in the matrix group, may offer clues as to how one measures distance in the complex plane. We note the partial derivatives in front of the trigonometric *sine* functions that act on *'r'*, represent changes in our complex radial component for the respective angles.

We can add any row or column together and the result will equal $e^{i(x+y+z)}$. We demonstrate below where we substitute $\begin{bmatrix} 0 & -1 \\ 1 & 0 \end{bmatrix} = i$ as $\mathbb{R}^2 \to \mathbb{C}^1$, where for the odd trigonometry functions

$$\left( \frac{2\cos(r)}{3} + \frac{\cos(\gamma)}{3} + \frac{\partial r}{\partial x} \cdot \frac{2i\sin(r)}{3} + \frac{i\sin(\gamma)}{3} - \frac{\cos(r)}{3} + \frac{\cos(\gamma)}{3} + \frac{\partial r}{\partial z} \cdot \frac{2i\sin(r)}{3} + \frac{i\sin(\gamma)}{3} - \frac{\cos(r)}{3} + \frac{\cos(\gamma)}{3} \frac{\partial r}{\partial y} \cdot \frac{2i\sin(r)}{3} + \frac{i\sin(\gamma)}{3} \right),$$

$$\to \cos(\gamma) + i\sin(\gamma) + \left( \frac{2\cos(r)}{3} - \frac{\cos(r)}{3} - \frac{\cos(r)}{3} \right) + \left( \frac{\partial r}{\partial x} \cdot \frac{2\sin(r)}{3} + \frac{\partial r}{\partial y} \cdot \frac{2\sin(r)}{3} + \frac{\partial r}{\partial z} \cdot \frac{2\sin(r)}{3} \right),$$

and as $\frac{\partial r}{\partial x} = \frac{(x - \frac{y}{2} - \frac{z}{2})}{r}$, $\frac{\partial r}{\partial y} = \frac{(-\frac{x}{2} + y - \frac{z}{2})}{r}$, and $\frac{\partial r}{\partial z} = \frac{(-\frac{x}{2} - \frac{y}{2} + z)}{r}$

$$\to \cos(\gamma) + i\sin(\gamma) = e^{i(x+y+z)}. \tag{17}$$

We recall that our matrix g *(13)* had zeroes along the main diagonal and with exponentiation, we have that our Lie group matrix *G* is "special" (S), as

$$\det(G) = 1 \tag{18}$$

and is orthogonal (O) as

$$G^T G = 1 \text{ and } G^T = G^{-1}. \tag{19}$$

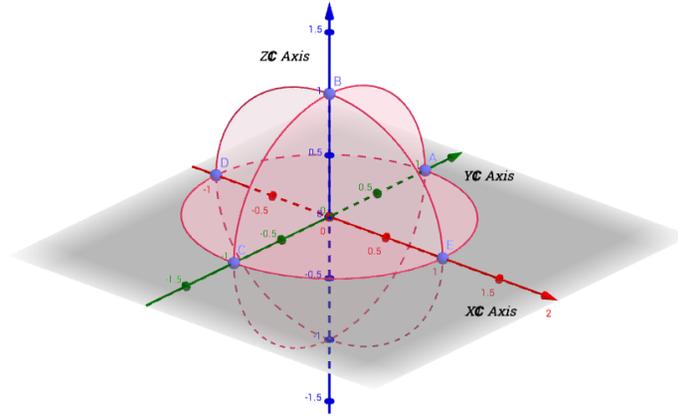

***Fig. 3*** – Graphic Representation of Axes of 6D Complex Unit Sphere

In ***Fig. 3***, we attempt to create a sense of the relationships between the angles on the complex unit sphere. We see it is an amalgam of the three $C_2$ cyclic groups that combined are $e^{ix+iy+iz}$. The three $C_2$ groups are arranged in orthogonal relation to one another and form a closed loop.

Calculating the eigenvalues of the matrix from ***Fig. 2***, we find

$$\lambda_1 = e^{-ir}, \lambda_2 = e^{-ir}, \lambda_3 = e^{ir}, \lambda_4 = e^{ir}, \lambda_5 = e^{-i\gamma}, \lambda_6 = e^{i\gamma} \tag{20}$$

Where $r = \sqrt{x^2 - xy - xz + y^2 - yz + z^2}$ and $\gamma = (x + y + z)$.

We see that the axes of rotation would be constantly changing based on the complex angle values we are calculating. This is similar to the eigenvalues one finds by calculating the same for the complex two-dimensional rotation matrix *(5)*. This implies an infinite number of rotation axes, and is much different from the eigenvalues the reader might be accustomed to in the form of fixed rotation axes.

Taking a step further, we can calculate our Lie group with the "c" and "e" permutations and we leave "a" out as this is identity permutation along the main diagonal that acts to uniformly scale the outputs. We fill in our g matrix from *(13)* to find

$$\mathfrak{g}_{(c,e)} = \begin{matrix} 0 & b & c & f & e & d \\ -b & 0 & -f & c & -d & e \\ e & f & 0 & d & c & b \\ -f & e & -d & 0 & -b & c \\ c & d & e & b & 0 & f \\ -d & c & -b & e & -f & 0 \end{matrix}. \tag{21}$$

We exponentiate the $\mathfrak{g}_{(c,e)}$ matrix above to find our Lie group. We demonstrate the first 2 orthonormal columns of the matrix below in ***Fig. 4***

$$\begin{bmatrix} \dfrac{e^{\left(-\frac{c}{2}-\frac{e}{2}\right)}\cdot 2\cos(r)}{3}+\dfrac{e^{(c+e)}\cdot\cos(\gamma)}{3} & \dfrac{\partial r}{\partial x}\cdot\dfrac{e^{\left(-\frac{c}{2}-\frac{e}{2}\right)}\cdot\sinh(r)}{3}+\dfrac{e^{(c+e)}\cdot\sin(\gamma)}{3} & \cdots \\ -\dfrac{\partial r}{\partial x}\cdot\dfrac{e^{\left(-\frac{c}{2}-\frac{e}{2}\right)}\cdot\sin(r)}{3}-\dfrac{e^{(c+e)}\cdot\sin(\gamma)}{3} & \dfrac{e^{\left(-\frac{c}{2}-\frac{e}{2}\right)}\cdot 2\cosh(r)}{3}+\dfrac{e^{(c+e)}\cdot\cos(\gamma)}{3} & \cdots \\ -\dfrac{e^{\left(-\frac{c}{2}-\frac{e}{2}\right)}\cdot\cos(r)}{3}+\dfrac{e^{(c+e)}\cdot\cos(\gamma)}{3} & \dfrac{\partial r}{\partial x}\cdot\dfrac{e^{\left(-\frac{c}{2}-\frac{e}{2}\right)}\cdot\sinh(r)}{3}+\dfrac{e^{(c+e)}\cdot\sin(\gamma)}{3} & \cdots \\ -\dfrac{\partial r}{\partial x}\cdot\dfrac{e^{\left(-\frac{c}{2}-\frac{e}{2}\right)}\cdot\sin(r)}{3}-\dfrac{e^{(c+e)}\cdot\sin(\gamma)}{3} & -\dfrac{e^{\left(-\frac{c}{2}-\frac{e}{2}\right)}\cdot\cosh(r)}{3}+\dfrac{e^{(c+e)}\cdot\cos(\gamma)}{3} & \cdots \\ -\dfrac{e^{\left(-\frac{c}{2}-\frac{e}{2}\right)}\cdot\cos(r)}{3}+\dfrac{e^{(c+e)}\cdot\cos(\gamma)}{3} & \dfrac{\partial r}{\partial x}\cdot\dfrac{e^{\left(-\frac{c}{2}-\frac{e}{2}\right)}\cdot\sinh(r)}{3}+\dfrac{e^{(c+e)}\cdot\sin(\gamma)}{3} & \cdots \\ -\dfrac{\partial r}{\partial x}\cdot\dfrac{e^{\left(-\frac{c}{2}-\frac{e}{2}\right)}\cdot\sin(r)}{3}-\dfrac{e^{(c+e)}\cdot\sin(\gamma)}{3} & -\dfrac{e^{\left(-\frac{c}{2}-\frac{e}{2}\right)}\cdot\cosh(r)}{3}+\dfrac{e^{(c+e)}\cdot\cos(\gamma)}{3} & \cdots \end{bmatrix}$$

***Fig. 4*** – First 2 Orthonormal Columns of Lie Group $G_{c,e}$.

In ***Fig. 5*** below we show an example of a waveform plot where we have set scaling factors for "c" and "e". We can see as the waveforms move from left to right, they are increasing in magnitude and becoming more symmetric. If we were to switch the signs of "c" and "e", we would see waveforms shrinking in magnitude as we progress from left to right. If "c" and "e" are equal in absolute value but of opposite sign, we would see no net change in the waveforms. Alternatively, if "c" and "e" are of the same sign their effects would add, and increase or decrease the waveform magnitude for positive or negative values respectively.

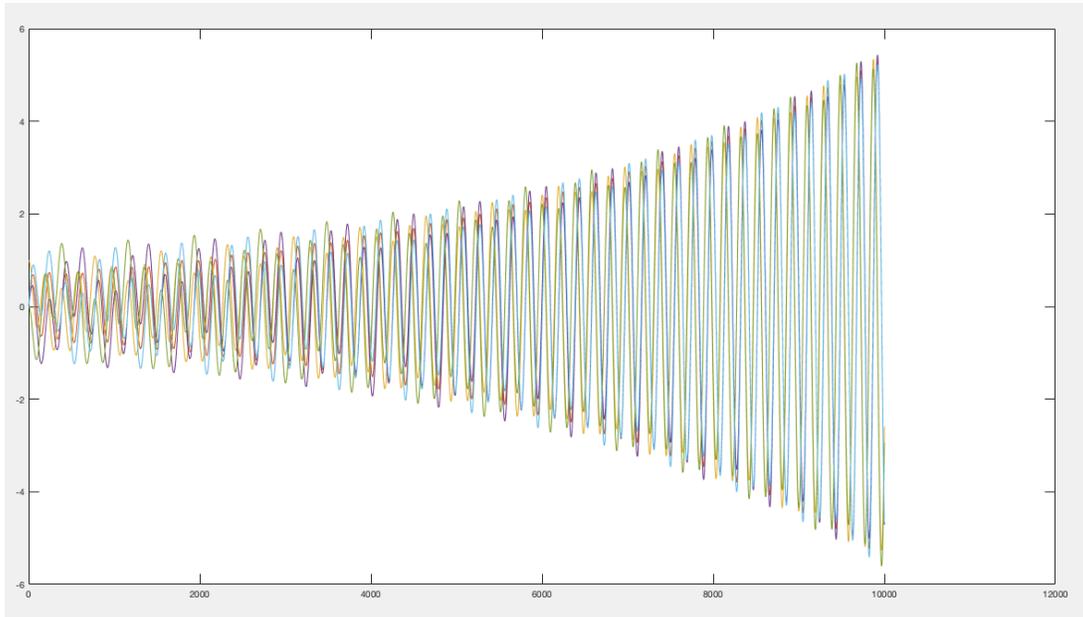

***Fig. 5*** – Waveform plot with scaling factors $c = \dfrac{-\pi}{30000}$ and $e = \dfrac{\pi}{10000}$

## Distance Functions

In regard to the Euclidean plane and that of space-time, physicists and mathematicians speak of distance functions that should equate to unity for "special" matrix groups. In the Euclidean case, we have become accustomed to taking the determinate of the two-dimensional complex rotation matrix and finding

$$\cos(x)^2 + \sin(x)^2 = 1. \tag{22}$$

In the hyperbolic rotation matrix, we find the determinant represents rotations in space-time otherwise known as space-time boost, and takes the form of [3]

$$\cosh(x)^2 - \sinh(x)^2 = 1. \tag{23}$$

By way of taking the determinate of a matrix group, we are able to find its distance function. Taking the distance function of the complex unit sphere's Lie group matrix $G$, we find a very nice reduction to our expected distance function for Euclidean space. To demonstrate, with taking the determinant we first find

$$\begin{aligned}\cos^2(\gamma)\cos^4(r) + \cos^2(\gamma)\sin^4(r) + \sin^2(\gamma)\cos^4(r) + \sin^2(\gamma)\sin^4(r) \\ + 2\cos^2(\gamma)\cos^2(r)\sin^2(r) + 2\sin^2(\gamma)\cos^2(r)\sin^2(r).\end{aligned} \tag{24}$$

This simplifies to

$$(\cos^2(\gamma) + \sin^2(\gamma))(\cos^2(r) + sin^2(r))^2. \tag{25}$$

Which ultimately reduces to

$$\begin{aligned}(\cos^2(\gamma) + \sin^2(\gamma)) &= 1 \\ \text{or}& \\ (\cos^2(r) + sin^2(r))^2 &= 1.\end{aligned} \tag{26}$$

To note, is that by taking the determinant for the Lie group matrix $G_{c,e}$, we find that it results in the same values as the above, where the "c" and "e" permutations cancel out in the calculation.

## The Unitary Representation U(3) of the Complex Unit Sphere

Taking what we have derived in our $\mathbb{R}^{2n}$ form of the six-dimensional Lie group we can transform this to the $\mathbb{C}^n$ representation and form our unitary group that is U(3). We can proceed keeping in mind that our prior imaginary permutations that acted as our group generators are representative of our respective imaginary entries as

$$\begin{bmatrix} 0 & 1 \\ -1 & 0 \end{bmatrix} = i. \tag{27}$$

We are in essence transforming $\mathbb{R}^{2n} \to \mathbb{C}^n$, or more specifically, $\mathbb{R}^6 \to \mathbb{C}^3$. We setup our permutations as such

$$b = \begin{bmatrix} i & 0 & 0 \\ 0 & 0 & i \\ 0 & i & 0 \end{bmatrix}, d = \begin{bmatrix} 0 & 0 & i \\ 0 & i & 0 \\ i & 0 & 0 \end{bmatrix}, f = \begin{bmatrix} 0 & i & 0 \\ i & 0 & 0 \\ 0 & 0 & i \end{bmatrix}. \tag{28}$$

Once again, we substitute x=b, y=d and z=f into the above, add the permutation together, and arrive at

$$g_{u(3)} = \begin{matrix} ix & iz & iy \\ iz & iy & ix \\ iy & ix & iz \end{matrix}. \tag{29}$$

To find our Lie group we can now simply exponentiate the matrix $\mathfrak{g}_{su(3)}$

$$G_{U(3)} = \exp(\mathfrak{g}_{u(3)}), \qquad (30)$$

which yields our U(3) Lie group

$$\begin{bmatrix} \frac{e^{i\gamma}}{3} + \frac{2\cos(r)}{3} + \frac{\partial r}{\partial x} \cdot \frac{2i\sin(r)}{3} & \frac{e^{i\gamma}}{3} - \frac{\cos(r)}{3} + \frac{\partial r}{\partial z} \cdot \frac{2i\sin(r)}{3} & \frac{e^{i\gamma}}{3} - \frac{\cos(r)}{3} + \frac{\partial r}{\partial y} \cdot \frac{2i\sin(r)}{3} \\ \frac{e^{i\gamma}}{3} - \frac{\cos(r)}{3} + \frac{\partial r}{\partial z} \cdot \frac{2i\sin(r)}{3} & \frac{e^{i\gamma}}{3} + \frac{2\cos(r)}{3} + \frac{\partial r}{\partial y} \cdot \frac{2i\sin(r)}{3} & \frac{e^{i\gamma}}{3} - \frac{\cos(r)}{3} + \frac{\partial r}{\partial x} \cdot \frac{2i\sin(r)}{3} \\ \frac{e^{i\gamma}}{3} - \frac{\cos(r)}{3} + \frac{\partial r}{\partial y} \cdot \frac{2i\sin(r)}{3} & \frac{e^{i\gamma}}{3} - \frac{\cos(r)}{3} + \frac{\partial r}{\partial x} \cdot \frac{2i\sin(r)}{3} & \frac{e^{i\gamma}}{3} + \frac{2\cos(r)}{3} + \frac{\partial r}{\partial z} \cdot \frac{2i\sin(r)}{3} \end{bmatrix}.$$

*Fig. 5* – $G_{U(3)}$ Lie group. $r = \sqrt{x^2 - xy - xz + y^2 - yz + z^2}$, $\gamma = (x + y + z)$, $\frac{\partial r}{\partial x} = \frac{(x - \frac{y}{2} - \frac{z}{2})}{r}$, $\frac{\partial r}{\partial y} = \frac{(-\frac{x}{2} + y - \frac{z}{2})}{r}$, $\frac{\partial r}{\partial z} = \frac{(-\frac{x}{2} - \frac{y}{2} + z)}{r}$.

We continue and calculate the determinant of $G_{U(3)}$ algebraically and find

$$\det(G_{U(3)}) = e^{i\gamma}\cos^2(r) + e^{i\gamma}\sin^2(r) = e^{i\gamma}(\cos^2(r) + \sin^2(r)) = e^{i\gamma}. \qquad (31)$$

The result above in *Eq. (31)* coincides with the fact that for any complex square matrix, we have that

$$\det(e^{\mathfrak{g}_{u(3)}}) = e^{Trace(\mathfrak{g}_{u(3)})}. \qquad (32)$$

We recall that $\gamma = (x + y + z)$, and we have that the modulus of the determinate will always be equal to one. This is a noteable result, as we normally require a trace equal to zero in our adjoint matrix to ensure it has a determinant equal to unity in the Lie group form. It is also easy to appreciate that the $G_{U(3)}$ Lie group is equal to its own transpose as $G_{U(3)} - G_{U(3)}^T = 0$ and is thus a symmetric matrix.

Many might decry that the SU(3) group used in physics and the Gell-Mann matrices are rather contrived by the lack of a division algebra derivation. We believe the above U(3) Lie group may bridge this gap and makes headway to prove the existence of the SU(3) Lie group mathematically and physically.

## Rotations of SO(3,C)

With regard to the dynamics of rotations of the complex unit sphere, we find the angles are coupled when doing spin rotations (complex Euclidean rotations). When you rotate by imaginary angles, the elements of the Lie group become hyperbolics. Rotations are oriented and as shown with the aforementioned verification of distance functions, we always arrive at unity for real-values angles. The rotation matrix is robust and we find we can rotate by single points, complex six-vectors and also complex six-dimensional volumes. Rotation angles can be real, imaginary as well as complex (both real and imaginary). Being non-abelian by multiplication, we find the arrangement of rotations does make a difference with the results we arrive at.

We find there are 2 main ways to rotate on the complex unit sphere. The first, is by rotation by a single angle $\angle X, \angle Y\ and\ \angle Z$, or by way of rotation by an oriented plane that is represented by "-c" and "-e" from our permutations.

We present the dynamics of rotations on the complex unit sphere, by first showing rotation by a single angle. We use the convention that a rotation by a positive real angle yields a real counter clock-wise rotation and our vector is multiplied to the right of our Lie group matrix. For these rotations about a single angle, we find the following relationships for a real unit volume being rotated about the respective angles $\angle X, \angle Y$ and $\angle Z$.

$$\begin{aligned}&\angle X \text{ Progression}\\ &X_r \to X_i \to -X_r \to -X_i \to X_r\\ &Y_r \to Z_i \to -Y_r \to -Z_i \to Y_r\\ &Z_r \to Y_i \to -Z_r \to -Y_i \to Z_r\end{aligned} \quad (32)$$

$$\begin{aligned}&\angle Y \text{ Progression}\\ &X_r \to Z_i \to -X_r \to -Z_i \to X_r\\ &Y_r \to Y_i \to -Y_r \to -Y_i \to Y_r\\ &Z_r \to X_i \to -Z_r \to -X_i \to Z_r\end{aligned} \quad (33)$$

$$\begin{aligned}&\angle Z \text{ Progression}\\ &X_r \to +Y_i \to -X_r \to -Y_i \to X_r\\ &Y_r \to +X_i \to -Y_r \to -X_i \to Y_r\\ &Z_r \to +Z_i \to -Z_r \to -Z_i \to Z_r\end{aligned} \quad (34)$$

Next, we present rotations of the oriented volumes, by the oriented planes of $-c$ and $-e$, will yield the following relations when an oriented volume $\hat{V}$ rotates by $\angle \frac{\pi}{2}$ show in **Tables 1.** and **2.**

### $-c$ oriented plane rotation progression

| Axis | $\hat{V}(-c^1)$ | $\hat{V}(-c^2)$ | $\hat{V}(-c^3)$ | $\hat{V}(-c^4)$ | $\hat{V}(-c^5)$ | $\hat{V}(-c^6)$ |
|---|---|---|---|---|---|---|
| $X \to$ | $-X_y$ | $X_z$ | $-X_x$ | $X_y$ | $-X_z$ | $X_x$ |
| $Y \to$ | $-Y_z$ | $Y_x$ | $-Y_y$ | $Y_z$ | $-Y_x$ | $Y_y$ |
| $Z \to$ | $-Z_x$ | $Z_y$ | $-Z_z$ | $Z_x$ | $-Z_y$ | $Z_z$ |

*Table 1.*

### $-e$ oriented plane rotation progression

| Axis | $\hat{V}(-e^1)$ | $\hat{V}(-e^2)$ | $\hat{V}(-e^3)$ | $\hat{V}(-e^4)$ | $\hat{V}(-e^5)$ | $\hat{V}(-e^6)$ |
|---|---|---|---|---|---|---|
| $X \to$ | $-X_z$ | $X_y$ | $-X_x$ | $X_z$ | $-X_y$ | $X_x$ |
| $Y \to$ | $-Y_x$ | $Y_z$ | $-Y_y$ | $Y_x$ | $-Y_z$ | $Y_y$ |
| $Z \to$ | $-Z_y$ | $Z_x$ | $-Z_z$ | $Z_y$ | $-Z_x$ | $Z_z$ |

*Table 2.*

We can notice that the above progressions will pass thru 6 different relationships before arriving back at the starting point when rotating by "$-c$" and "$-e$". The "$-c$" rotation will proceed in an overall clockwise fashion, but it is important to note there is alternation between positive and negative axes. So, for "$-c$", we see that X,Y,Z represents the order or successive rotations, whereas with "$-e$", we find the order opposite and ordered as Z,Y,X.

One very profound relationship that deeper inspection reveals, is that when we are rotating on the complex unit sphere from $0$ to $2\pi$, we have the same result of net-zero distance when a complete revolution of $2\pi$ is made. This is related to the same phenomena we see on the complex unit circle where

$$\int_0^{2\pi} e^{ix} dx = 0. \quad (35)$$

The reader can see for both "-c" and "-e" plane rotations, starting with a unit volume that is all real, it will become negative and opposite to the initial starting for $\hat{V}(-e^3)$ and $\hat{V}(-c^3)$ rotations. The coresponding integral for *Eqn. (35)* from above in terms of the complex unit sphere is

$$\iiint_0^{2\pi} e^{ix+iy+iz}\, dx\, dy\, dz = 0. \tag{36}$$

**Plotting on the Complex Unit Sphere**

We continue by looking into some basic plots on the complex unit sphere. We define our six-vector $\hat{v}$ as

$$\hat{v}^T = (1,0,0,0,0,0). \tag{37}$$

We first confirm that the SO(3,C) form reduces to U(1) by rotating $\hat{v}$ by the $\angle X = \frac{\pi}{1000}$ and iterate five thousand times. We see in *Fig. 6*, a perfect circular form and the familiar cosine and sine waveforms.

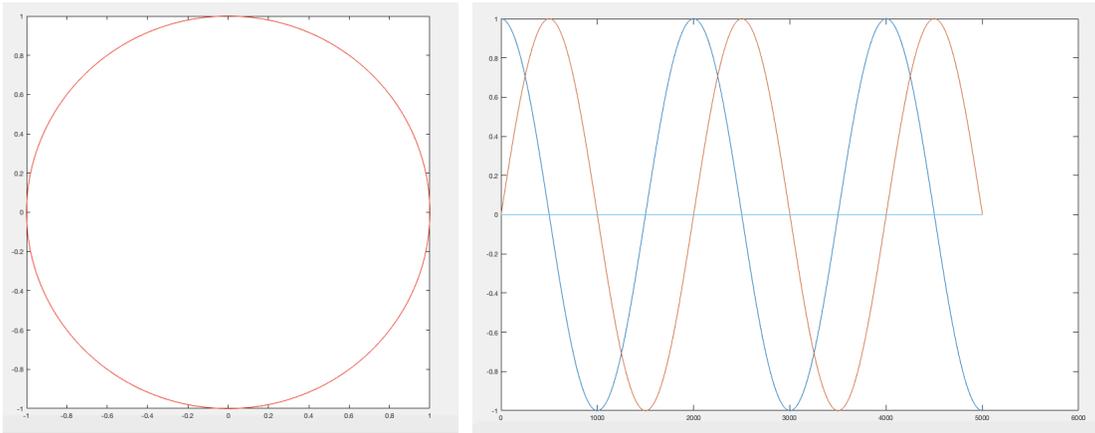

*Fig. 6* – Geodesic on the complex unit sphere – plot of rotation by angle $\angle X$ which reduces to $e^{ix}$ and its waveforms, $\cos(x)$ and $i\sin(x)$.

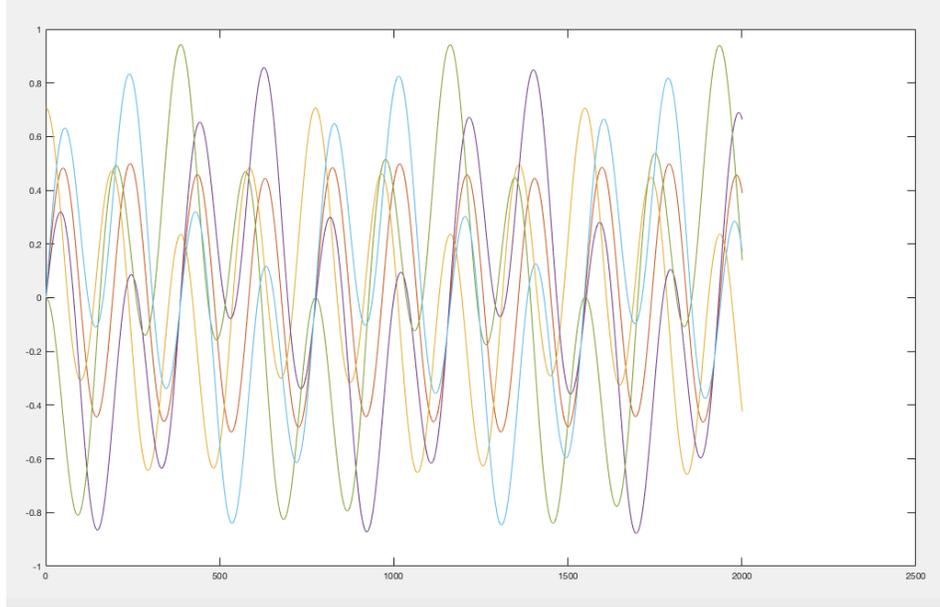

***Fig. 7*** – Composite rotation waveforms – rotation angles $\angle X = \frac{\pi}{200}, \angle Y = \frac{\pi}{300}, \angle Z = \frac{\pi}{500}$

In ***Fig. 7*** above we rotate our vector $\hat{v}$ by $\angle X = \frac{\pi}{200}, \angle Y = \frac{\pi}{300}, \angle Z = \frac{\pi}{500}$, and iterate two thousand times. We can see the inter-relation amongst all six waveforms over time.

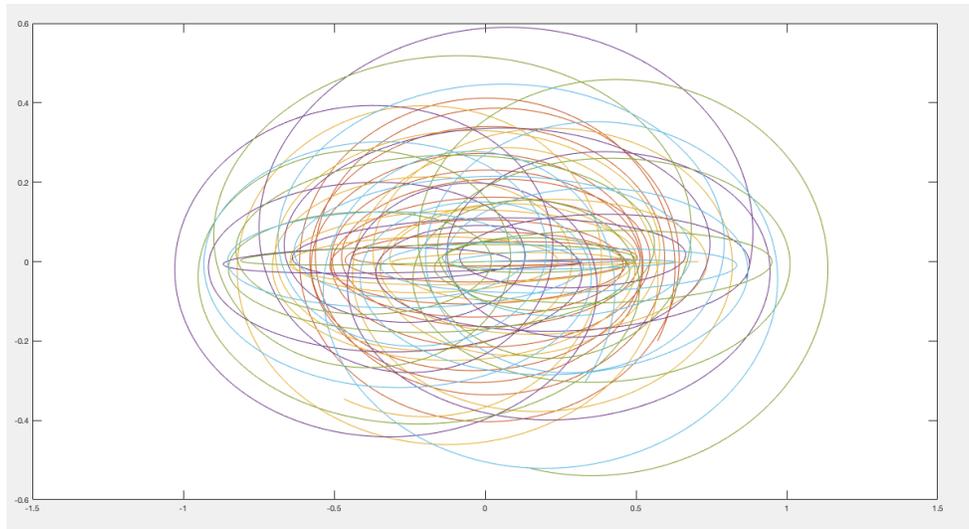

***Fig. 8*** – Complex Composite rotation waveforms– rotation angles $\angle X = \frac{\pi}{200} + \frac{i\pi}{10000}, \angle Y = \frac{\pi}{300} + \frac{i\pi}{20000}, \angle Z = \frac{\pi}{500} + \frac{i\pi}{30000}$.

In ***Fig. 8*** we rotate our vector $\hat{v}$ by the complex and composite angles $\angle X = \frac{\pi}{200} + \frac{i\pi}{10000}, \angle Y = \frac{\pi}{300} + \frac{i\pi}{20000}, \angle Z = \frac{\pi}{500} + \frac{i\pi}{30000}$ thru two thousand iterations. We can see the waveforms that we noted before, start to become elliptical from the effects of the imaginary angles that we added.

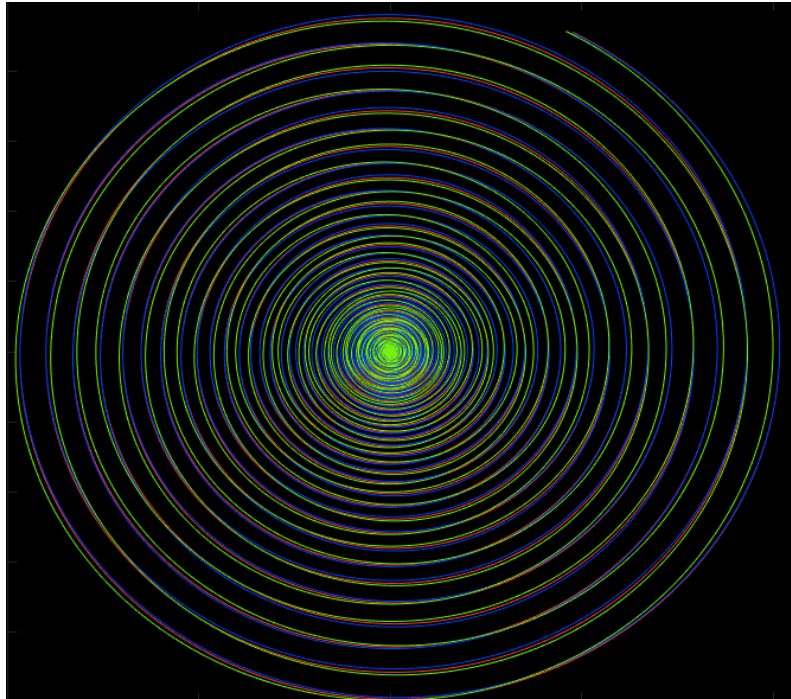

***Fig. 9*** – Complex Composite rotation waveforms– rotation angles $\angle X = \frac{\pi}{200} + \frac{i\pi}{10000}, \angle Y = \frac{\pi}{300} + \frac{i\pi}{20000}, \angle Z = \frac{\pi}{500} + \frac{i\pi}{30000}$ Twenty-thousand iterations.

In ***Fig. 9*** we continue the same rotation angles from ***Fig. 8*** and iterate for twenty-thousand times to see that the waveforms have become more elliptical as we move away from the origin.

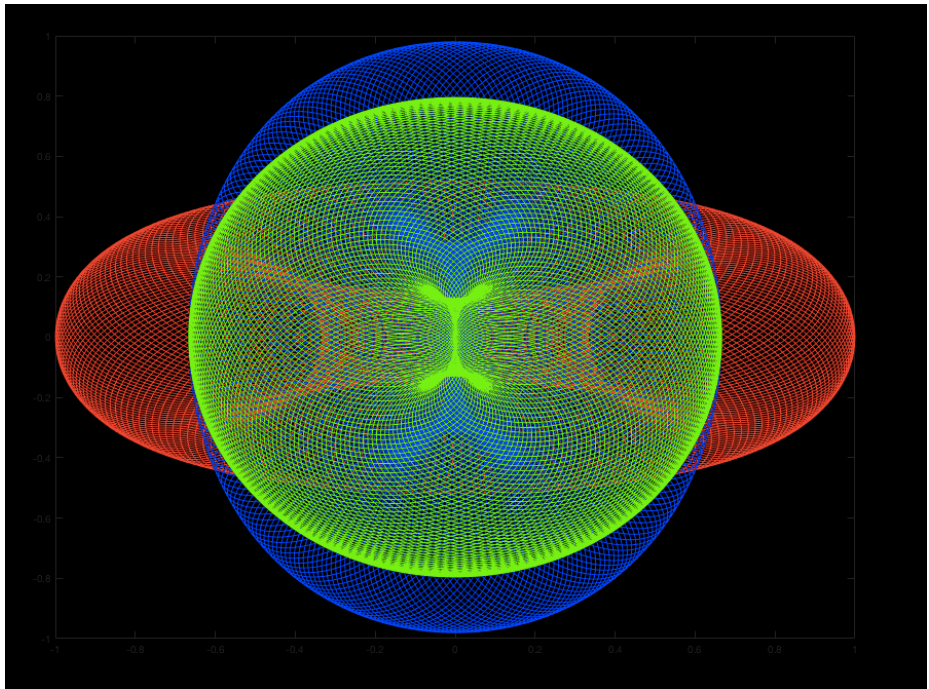

***Fig. 10*** – Composite rotation plot – Toroidal forms $\angle X = \frac{\pi}{200}, \angle Y = \frac{\pi}{300}, \angle Z = \frac{\pi}{500}$.

In **_Fig. 10_**, we rotate our vector $\hat{v}$ by the angles $\angle X = \frac{\pi}{200}, \angle Y = \frac{\pi}{300}, \angle Z = \frac{\pi}{500}$ over twenty-five thousand iterations. We plot each respective 2-D form that is X, Y and Z. We notice beautiful forms that take shape in the guise of toroidal shapes for each.

**Conclusion**

By using a very simple approach that started with using selective permutation matrices of a finite group and by way of exponentiation, we are able to form the Lie groups of the complex unit sphere SO(3,C) and reduce this to its unitary form in U(3). Rotating with these Lie groups is straightforward, robust and takes advantage of exponential functions. With these higher dimensional Lie groups, we are able to find our familiar distance functions, demonstrate complex angular rotations and all the while being reassured that the derivations represent division algebras. We can entertain the possibility that these Lie groups represent how rotations may occur in empty space, that each plane of rotation requires a respective $C_2$ algebra, and to make rotations in what we perceive as three-dimensions, requires the orthogonal arrangement of three $C_2$ algebras to do so. We believe the aforementioned derivations have the potential to advance our approaches to relativity, quantum field theory, complex analysis, deriving higher dimensional forms of similar division algebras and technologies requiring orientated rotations.